\documentclass[12pt]{article}

\usepackage{graphicx,tikz,color,url}
\usepackage{amsmath,amssymb,latexsym}
\usepackage{pgfplots}
\usepackage{mathrsfs}
\usetikzlibrary{arrows}
\pgfplotsset{compat=1.15}

\newtheorem{theorem}{Theorem}[section]

\newtheorem{lemma}[theorem]{Lemma}
\newtheorem{proposition}[theorem]{Proposition}
\newtheorem{observation}[theorem]{Observation}
\newtheorem{conjecture}[theorem]{Conjecture}

\newcommand{\proof}{\noindent{\bf Proof.\ }}
\newcommand{\qed}{\hfill $\square$ \bigskip}

\newcommand{\gsmb}{\gamma_{\rm SMB}}
\newcommand{\gmb}{\gamma_{\rm MB}}

\newcommand{\cH}{{\cal H}}

\newcommand{\cD}{{\cal D}}

\newcommand{\w}{{\rm w}}

\begin{document}

\title{Fast winning strategies for Staller in the Maker-Breaker domination game}

\author{Csilla Bujt\'as\thanks{Email: \texttt{csilla.bujtas@fmf.uni-lj.si}}  \qquad \qquad
	Pakanun Dokyeesun\thanks{Email: \texttt{pakanun.dokyeesun@student.fmf.uni-lj.si}}
}
\date{\empty}
\maketitle
\vspace{-15pt}

\begin{center}
Faculty of Mathematics and Physics\\ University of Ljubljana, Slovenia\\
	\medskip
	
\end{center}	

\begin{abstract}
	The Maker-Breaker domination game is played on a graph $G$ by two players, called Dominator and Staller, who alternately choose a vertex that has not been played so far. Dominator wins the game  if his moves form a dominating set. Staller wins if she plays all vertices from a closed neighborhood of a vertex $v \in V(G)$. Dominator's fast winning strategies were studied earlier. In this work, we concentrate on the cases when Staller has a winning strategy in the game. We introduce the invariant $\gamma'_{\rm SMB}(G)$ (resp., $\gamma_{\rm SMB}(G)$) which is the smallest integer $k$ such that, under any strategy of Dominator, Staller can win the game by playing at most $k$ vertices, if Staller (resp., Dominator) plays first on the graph $G$.
	
	We prove some basic properties of $\gamma_{\rm SMB}(G)$ and $\gamma'_{\rm SMB}(G)$ and study the parameters' changes under some operators as taking the disjoint union of graphs or deleting a cut vertex. We show that the inequality $\delta(G)+1 \le \gamma'_{\rm SMB}(G) \le \gamma_{\rm SMB}(G)$ always holds and that for every three integers $r,s,t$ with $2\le r\le s\le t$, there exists a graph $G$ such that $\delta(G)+1 = r$, $\gamma'_{\rm SMB}(G) = s$, and $\gamma_{\rm SMB}(G) = t$. We prove exact formulas for $\gamma'_{\rm SMB}(G)$ where $G$ is a path, or it is a tadpole graph which is obtained from the disjoint union of a cycle and a path by adding one edge between them. 
\end{abstract}
\noindent
{\bf Keywords:} domination game; Maker-Breaker game; winning number; Maker-Breaker domination game; closed neighborhood hypergraph\\


\section{Introduction}
\label{sec:intro}

A general frame for studying Maker-Breaker games was introduced by Erd\H os and Selfridge \cite{erdos-1973} nearly $50$ years ago. Since then, several variants have been defined and studied extensively \cite{hefetz}. In this work, we are interested in the Maker-Breaker domination game \cite{duchene-2019+}. In particular, we focus on the cases when Breaker (i.e., Staller) wins the game and consider Breaker's fast winning strategies. Along this way, we translate the problem to the Maker-Breaker game played on the closed neighborhood hypergraph of a graph. We establish sharp lower and upper bounds on the corresponding graph invariants and determine the exact values for paths and tadpole graphs.

\subsection{Standard definitions}

Throughout the paper, we denote by $V(G)$ and $E(G)$, respectively, the vertex and edge set of a graph $G$. For a vertex $v \in V(G)$, the \emph{open neighborhood}  $N_G(v)$ contains the neighbors of $v$ and the \emph{closed neighborhood} of $v$ is defined as $N_G[v]=N_G(v) \cup \{v\}$. Then, $\deg_G(v)=|N_G(v)|$. For a set $S \subseteq V(G)$ we use the analogous notation $N_G(S)= \bigcup_{v \in S}N_G(v)$ and $N_G[S]= \bigcup_{v \in S}N_G[v]$. 
If the graph is clear from the context, we omit the index and simply write $N(v)$, $N[v]$, $N(S)$, and $N[S]$. A vertex $v$ is a \textit{leaf} if $\deg_G(v)=1$, and the neighbor of a leaf is a  \textit{support vertex}. A support vertex is said to be a \textit{strong support vertex} if it is adjacent to at least two leaves. Otherwise, we call it a \textit{weak support vertex}. 

  A vertex dominates itself and its neighbors.  A set $D \subseteq V(G)$ is a \emph{dominating set} in $G$ if $N[D]=V(G)$. The minimum cardinality of such a set is the \emph{domination number} of $G$ and denoted by $\gamma(G)$. A dominating set is \emph{minimal} if it does not contain any dominating set as a proper subset. 

A hypergraph $\cH$ is a set system over the underlying vertex set $V(\cH)$. Equivalently, the set $E(\cH)$ of its (hyper)edges is a subset of $2^{V(\cH)}$. In general, we assume that each $e \in E(\cH)$ contains at least one vertex, but sometimes we do not exclude the case of $e=\emptyset$. For the latter situation, we will always mention this possibility directly. A hypergraph is \emph{$k$-uniform} if $|e|=k$ for every $e \in E(\cH)$ and it is \emph{simple} if $e \not\subseteq f$ holds for every two different edges. Note that a simple graph is just a simple $2$-uniform hypergraph. We say that  $\cH_1$ is a \emph{subhypergraph} of $\cH_2$, if $E(\cH_1) \subseteq E(\cH_2)$.

A vertex set $T \subseteq V(\cH)$ is a \emph{transversal} (or vertex cover) in a hypergraph $\cH$ if it meets every edge of $\cH$. Assuming $\emptyset \notin E(\cH)$, the minimum cardinality of a transversal is the \emph{transversal number} $\tau(\cH)$ of the hypergraph. A \emph{minimal transversal} of $\cH$ is a transversal that does not contain any other transversal as a proper subset. Following Berge's notation, we denote by $Tr(\cH)$ the hypergraph on the vertex set $V(\cH)$ that contains all minimal transversals of $\cH$ as hyperedges. We cite a related result of Berge.
\begin{proposition}[Corollary~1, pp.\ 44 \cite{berge}] \label{prop:berge}
Let $\cH_1$ and $\cH_2$ be two simple hypergraphs on the same vertex set. Then, $Tr(\cH_1)=\cH_2$ if and only if $Tr(\cH_2)=\cH_1$.
\end{proposition}

We also introduce an operator, named \emph{simplification},  on hypergraphs that transforms $\cH$ into $\widehat{\cH}$ by sequentially removing, while it is possible, an edge $e \in E(H)$ if it contains another edge as a subset. It is clear by definition that  $\widehat{\cH}$ is a simple hypergraph and that $T$ is a (minimal) transversal in $\cH$ if and only if it is a (minimal) transversal in $\widehat{\cH}$.

Given a graph $G$, we define two hypergraphs on the same vertex set $V(G)$.  Let $\cD_G$ be the hypergraph that contains all the minimal dominating sets of $G$ as hyperedges and let $\cH_G$ be the hypergraph with edge set $E(\cH_G)=\{N_G[v]: v \in V(G)\}$. According to  the standard terminology, we say that $\cH_G$ is the  \emph{closed neighborhood hypergraph} of $G$. By definition, $D \subseteq V(G)$ is a minimal dominating set in $G$ if and only if it is a minimal transversal in the closed neighborhood hypergraph $\cH_G$. Observe that $\cD_G$ is a simple hypergraph for every $G$. These observations together with Proposition~\ref{prop:berge} imply the following statement.
\begin{proposition} \label{prop:H-and-D}
It holds for every graph $G$ that $Tr(\cH_G)=\cD_G$ and $Tr(\cD_G)=\widehat{\cH_G}$.
\end{proposition}

\subsection{Winning number in a Maker-Breaker game}

The \textit{Maker-Breaker game} is played by two persons, namely Maker and Breaker, on a hypergraph $\cH$. The vertex set $V(\cH)$ is called the \emph{board} of the game while the edges of $\cH$ are interpreted as  \textit{winning sets}. The two players alternately choose (i.e., play) an unplayed vertex from the board. We say that a player \emph{claims} the set $S \subseteq V(\cH)$ if he (or she) plays all vertices of $S$. It is a \emph{Maker-start game} (resp., \emph{Breaker-start game}) if Maker (resp., Breaker) is the first to play.
\emph{Maker wins} the game if he claims a winning set while \emph{Breaker wins} if she can prevent Maker from doing this. The latter equivalently means that Breaker plays at least one vertex from each winning set, that is, she claims a minimal transversal of $\cH$. We may conclude the following.
\begin{observation} \label{obs:switch}
   A Maker-Breaker game on $\cH$ with a player A as Maker and a player B as Breaker is the same as the Maker-Breaker game on $Tr(\cH)$ where the roles of A and B are switched.
 \end{observation}

In the last decades, several variants of the Maker-Breaker game were specified and intensively studied. These include the connectivity game \cite{gebauer-2010, hefetz-2012}, the hamiltonicity game \cite{ben-2011,Krivelevich-2010,stojakovic}, the Maker-Breaker resolving \cite{kang-2021+}, domination \cite{duchene-2019+, gledel-2019} and total domination \cite{gledel-2019+} games. For other versions and references see \cite{hefetz}.

\medskip

Fast winning strategies were also studied for some types of Maker-Breaker games \cite{clem-2012, clem-2018,duchene-2019+,gledel-2019+,gledel-2019,kang-2021+}. We introduce here four hypergraph invariants to measure how fast a player can win in a (general) Maker-Breaker game. Assume first that a Maker-Breaker game is played on $\cH$ and Maker can win the game when he is the first player. We say that Maker wins in his $i^{\rm th}$ move if he claims a winning set (i.e., an edge of $\cH$) with this move. We suppose that Maker's goal is to win the game as soon as possible and that Breaker's goal is just the opposite. The \emph{winning number of Maker}, denoted by $\w_M^M(\cH)$, is the minimum number of his moves he needs to win the game if both players play optimally. If Maker has no winning strategy as a first player, we set $\w_M^M(\cH)=\infty$. When it is a Breaker-start game, $\w_M^B(\cH)$ is defined similarly.  

Now, consider the problem from Breaker's point of view. We say that the game is over and Breaker wins in her $i^{\rm th}$ move if she claims a minimal transversal of $\cH$ with this move. Clearly, this situation makes it impossible for Maker to claim a winning set even if the game would be continued. Analogously to the previous case, we assume that Breaker's goal is to win the game by playing as few vertices as possible, while Maker's goal is the opposite. Then, if Maker starts the game and Breaker has a winning strategy, the \emph{winning number of Breaker},  $\w_B^M(\cH)$, is the minimum number of her moves needed to win the game if both players play optimally. Further,  $\w_B^M(\cH)=\infty$, if Breaker cannot win the game.  For the game when Breaker is the a first player,  $\w_B^B(\cH)$ is defined analogously. Concerning the notation $\w_X^Y$ in general, for every pair $(X,Y)$ with $X,Y \in \{M,B\}$ the superscript $Y$ denotes the player who starts the game while the subscript $X$ indicates the player whose fast winning strategies are considered.  

If $\w_X^Y(\cH)$ equals $\infty$ (for $X,Y \in \{M,B\}$) by definition, we consider it as an infinite cardinality number (e.g., $\aleph_0$) and refer to it in calculations or relations accordingly. For instance, for every hypergraph $\cH$ and $X\in \{M,B\}$,  $\w_M^X(\cH)+\w_B^X(\cH)=\infty$ always holds,  while $\w_M^X(\cH) < \w_B^X(\cH)$ means that Maker has a winning strategy in the game if $X$ plays first. 
By Observation~\ref{obs:switch}, $\w_M^M(\cH)=\w_B^B(Tr(\cH))$ and $\w_M^B(\cH)=\w_B^M(Tr(\cH))$.

\subsection{Maker-Breaker domination game}
The \textit{Maker-Breaker domination game} (\emph{MBD game} for short) was recently introduced in \cite{duchene-2019+} and further studied in \cite{gledel-2019}. The MBD game on a graph $G$ is a Maker-Breaker game played on $\cD_G$. Following the terminology introduced for the earlier versions of the domination game \cite{borowiecki-2019, all-2019, bresar-2010, book-2021, bujtas-2016, henning-2015, henning-2017a, kinnersley-2013}, Maker and Breaker are respectively called Dominator and Staller in an MBD game. 
An MBD game is called \emph{D-game} if Dominator starts the game and it is an \emph{S-game} if Staller starts. 
By Observation~\ref{obs:switch} and Proposition~\ref{prop:H-and-D}, a Maker-Breaker domination game on a graph $G$ can also be considered as a Maker-Breaker game on the closed neighborhood hypergraph $\cH_G$ (or on its simplification $\widehat{\cH_G}$) where Staller is the Maker and Dominator is the Breaker.

Now suppose that Dominator has a winning strategy on $G$ in the D-game. As it was introduced in \cite{gledel-2019}, the minimum number of moves of Dominator to win the game is the \textit{Maker-Breaker domination number} (\emph{MBD-number}) $\gmb(G)$. If Dominator does not have a winning strategy, we set $\gmb(G) = \infty$. In terms of the winning number defined for the general case, $\gmb(G)=\w_M^M(\cD_G)$. Similarly, for the S-game,  $\gmb'(G)=\w_M^B(\cD_G)$.

In this paper, as proposed in \cite{gledel-2019}, we consider the problem from Staller's point of view. Suppose that Staller has a strategy to win the (MBD) D-game on a graph $G$. Then, Staller wins in her $i^{\rm th}$ move if she claims a minimal transversal of $\cD_G$ with this move. It equivalently means that she can claim a closed neighborhood $N_G[v]$ of a vertex $v$. The \emph{Staller-Maker-Breaker domination number} (\emph{SMBD-number})  $\gsmb(G)$ of $G$ is the minimum number of Staller's moves she needs to win the D-game if both players play optimally. If Staller has no winning strategy in the D-game, we set $\gsmb(G)=\infty$. This definition gives $\gsmb(G)= \w_B^M(\cD_G)=\w_M^B(\cH_G)$. For the S-game, the invariant $\gsmb'(G)$ is defined analogously and we have $\gsmb'(G)= \w_B^B(\cD_G)=\w_M^M(\cH_G)$.

By definition, a Maker-Breaker game is a finite game and there is always a winner. Therefore, for each graph $G$ either both $\gmb(G) < \infty$ and $\gsmb(G)=\infty$ hold, or $\gmb(G)=\infty$ and $\gsmb(G) <\infty$. The analogous statement is true for the invariants related to the S-game. Further, as it was shown in \cite[Proposition~2.1.6, pp.\ 15]{hefetz-2012} for the Maker-Breaker game in general, if Staller (Dominator) can win as the second player, then she (he) can win as the first player. Then, $\gsmb(G) < \infty$ implies $\gsmb'(G) < \infty$. 

In this work, we concentrate on the cases where Staller wins in at least one of the D- and S-game that is if $\gsmb'(G) < \infty$. Typically, we will consider the Maker-Breaker game on the closed neighborhood hypergraph $\cH_G$.
\bigskip

\paragraph{Structure of the paper.} In Section~\ref{sec:prelim-1} we study the behavior of the winning numbers under some changes in the hypergraph such as taking disjoint union, deleting or shrinking edges. Then, in Section~\ref{sec:prelim-2}, we state similar results for the SMBD-numbers of a graph and prove the inequality $\delta(G) +1 \le \gsmb'(G) \le \gsmb(G)$ and the No-Skip Lemma. Section~\ref{sec:relation} is devoted to a construction which shows that for every three integers satisfying $2 \le r \le s \le t$ there exists a graph $G$ with $\delta(G)+1 = r$, $\gsmb(G) = s$, and $\gsmb'(G) = t$. 
In Section~\ref{sec:support}, we consider the changes of SMBD-numbers if a leaf and its support vertex is removed from a graph. Applying the previous propositions, it is proved  in Section~\ref{sec:path} that $\gsmb'(P_n)=\lfloor \log_2 n \rfloor +1 $ holds for every path of odd order. The SMBD-numbers of tadpole graphs (that can be obtained from a cycle and a path by adding an edge between a vertex of the cycle and an endvertex of the path) are determined in Section~\ref{sec:tadpole}. Some additional observations and conjectures are provided in Section~\ref{sec:concluding}.

\section{Basic properties}
\label{sec:prelim}

In this section we state several basic properties of the Maker-Breaker games in general, on the MBD games, and their hypergraph representation. 

\subsection{Basics for Maker-Breaker games} \label{sec:prelim-1}
We introduce two operators on hypergraphs that help us modeling Maker's and Breaker's moves. Let $\cH$ be a hypergraph and $X \subseteq V(\cH)$. The hypergraph $\cH-X$ is obtained from $\cH$ by deleting the vertices contained in $X$ and also deleting all incident edges. Formally, 
$$V(\cH-X)=V(\cH)\setminus X \qquad \mbox{and} \qquad  
E(\cH-X)=\{e: e \in E(\cH) \enskip  {\rm and }\enskip e \cap X=\emptyset\}.$$
The other operator is named \emph{shrinking} and the resulted hypergraph is denoted by $\cH\mid X$. It means that the vertices in $X$ are deleted from $\cH$ but we keep the remaining parts of the incident edges. That is, 
$$V(\cH\mid X)=V(\cH)\setminus X \qquad \mbox{and} \qquad
 E(\cH\mid X)=\{e \setminus X: e \in E(\cH)\}.$$
Note that shrinking may transform a simple hypergraph into a non-simple one. In particular, it may create some empty edges. 
If $X$ contains only one vertex $v$, we may write $\cH-v$ and $\cH \mid v$ instead of $\cH- \{v\}$ and $\cH\mid \{v\}$, respectively.

\medskip
In a Maker-Breaker game, the winning sets can be modified after each move of the players as follows.
\begin{proposition} \label{prop:delete-shrink}
	Suppose that a Maker-Breaker game is played on a hypergraph $\cH$ and $v \in V(\cH)$.
	\begin{itemize}
	 \item[$(i)$] If Breaker is the first player and she plays $v$, the continuation of the game corresponds to a Maker-start game on $\cH-v$. In particular, if $v$ is an optimal first move of Breaker, then $\w_M^B(\cH)= \w_M^M(\cH -v)$.
	 \item[$(ii)$] If Maker is the first player and he plays $v$, the continuation of the game corresponds to a Breaker-start game on $\cH \mid v$. In particular, if $v$ is an optimal first move of Maker, then $\w_M^M(\cH)= \w_M^B(\cH \mid v)+1$. Further, Maker wins the game with this move $v$ if and only if $\cH \mid v$ contains an empty hyperedge but $\cH$ does not.
	 
	\end{itemize}
\end{proposition}
\proof $(i)$ After Breaker plays $v$, this vertex cannot be played anymore, and none of the incident edges can be claimed by Maker. As the other winning sets remain untouched, Maker may imagine that he starts a Maker-Breaker game on $\cH-v$ with his next move. Concerning the equality $\w_M^B(\cH)= \w_M^M(\cH -v)$, recall that only Maker's moves are counted for $\w_M^B(\cH)$.

$(ii)$ If Maker plays $v$ on $\cH$, this vertex cannot be played in the continuation. Further, if $v$ belongs to a winning set $e \in E(\cH)$ then, in the continuation, it is enough to claim $e\setminus \{v\}$ for Maker to win the game. Thus, the players may imagine that the continuation is a Breaker-start game on $\cH \mid v$.  It is also clear that Maker claims a winning set with the move $v$ if and only if $\{v\} \in E(\cH)$ that is equivalent to $\emptyset \in E(\cH \mid v)$.
\qed

Next we show that adding some new winning sets or shrinking some of the winning sets is never disadvantageous for Maker.
\begin{proposition}
	\label{prop:delete}
Let $\cH_1$ and $\cH_2$ be two hypergraphs on the same vertex set. 
\begin{itemize}
\item[$(i)$] If $E(\cH_1) \subseteq E(\cH_2)$, then $\w_M^M(\cH_1)\ge  \w_M^M(\cH_2)$ and $\w_M^B(\cH_1)\ge  \w_M^B(\cH_2)$ holds.
\item[$(ii)$] Suppose that for each $e \in E(\cH_1)$ there exists an edge $e'\in E(\cH_2)$ such that $e'\subseteq e$. Then, $\w_M^M(\cH_1)\ge  \w_M^M(\cH_2)$ and $\w_M^B(\cH_1)\ge  \w_M^B(\cH_2)$ holds.
\end{itemize}
\end{proposition}
\proof As $(i)$ can be considered as a special case of $(ii)$, it is enough to prove the second statement. If $\w_M^M(\cH_1)=\infty$ or $\w_M^B(\cH_1)=\infty$, the corresponding inequality clearly holds. Otherwise, suppose that Maker plays on $\cH_2$ (as the first player) following a strategy that is optimal for $\cH_1$ while Breaker plays optimally on $\cH_2$. Maker's strategy ensures that, under any strategy of Breaker, he claims a winning set from $E(\cH_1)$ in at most $\w_M^M(\cH_1)$ moves. Let this winning set be $e\in E(\cH_1)$ in the current game. Since there is a winning set $e'$ in $\cH_2$ so that $e' \subseteq e$, Maker wins on $\cH_2$ either with this move or earlier. This proves $\w_M^M(\cH_1)\ge  \w_M^M(\cH_2)$. The argument is similar if Breaker starts the game. \qed

Propositions~\ref{prop:delete-shrink} and \ref{prop:delete} together imply that, during a Maker-Breaker game, skipping a move or playing an isolated vertex (i.e., a vertex of degree $0$ in $\cH$) is never advantageous for Maker. The same is true for Breaker concerning her moves.

\begin{proposition}
	\label{prop:H-comp}
	Let $\cH$ be a disconnected hypergraph that consists of the components $\cH_1, \dots, \cH_\ell$ such that  $\w_M^M(\cH_1) \le \dots \le \w_M^M(\cH_\ell)$. Then, the following holds:
$$\w_M^M(\cH) =\w_M^M(\cH_1) \qquad  {\rm and} \qquad \w_M^M(\cH_1) \le \w_M^B(\cH) \le \w_M^M(\cH_2).$$
\end{proposition}		

\proof Since $E(\cH_1) \subseteq E(\cH)$, Proposition~\ref{prop:delete} (i) implies $\w_M^M(\cH) \le \w_M^M(\cH_1)$. On the other hand, Breaker may reply to each move of Maker by playing an optimal move in the same component. This strategy ensures the game on $\cH$ not being finished before Maker's $\w_M^M(\cH_1)^{\rm th}$ move. We infer of $\w_M^M(\cH) =\w_M^M(\cH_1)$.

By Proposition~\ref{prop:delete-shrink} (i) and Proposition~\ref{prop:delete} (i), $\w_M^M(\cH) \le \w_M^B(\cH)$ always holds. Thus, $\w_M^M(\cH) =\w_M^M(\cH_1)$ implies $\w_M^M(\cH_1) \le \w_M^B(\cH)$. For the upper bound, if Breaker starts the game on $\cH$ by selecting a vertex $v$,  consider the following strategy of Maker. If $v\in V(\cH_1)$,  Maker plays according to an optimal strategy on $\cH_2$ from the next move on. This ensures that the game finishes in his $\w_M^M(\cH_2)^{\rm th}$ move or earlier. If $v\notin V(\cH_1)$ then, from the next move, Maker plays according to an optimal strategy on $\cH_1$ and ensures that the game finishes in his $\w_M^M(\cH_1)^{\rm th}$ move or earlier. As $\w_M^M(\cH_1) \le \w_M^M(\cH_2)$, we may conclude $\w_M^B(\cH) \le \w_M^M(\cH_2)$. \qed

\subsection{Basics for Maker-Breaker domination games} \label{sec:prelim-2}

Recall that an MBD game on a graph $G$ corresponds to a Maker-Breaker game on the closed neighborhood hypergraph $\cH_G$ with $\gsmb'(G)=\w_M^M(\cH_G)$ and $\gsmb(G)=\w_M^B(\cH_G)$. When we consider an MBD game in which Staller has already played the vertices from $S$ and Dominator has played the vertices from $D$, we will often refer to the hypergraph of current winning sets. By Proposition~\ref{prop:delete-shrink}, this hypergraph is $(\cH_G \mid S)-D$. To prove lower bounds on $\gsmb(G)$ or $\gsmb'(G)$ we often modify the hypergraph $\cH$ of current winning sets by adding new winning sets or replacing an $e\in E(\cH)$ with an $e' \subseteq e$. By Proposition~\ref{prop:delete}, it makes easier for Staller to win. In the other case,  when we want to prove an upper bound on $\gsmb(G)$ or $\gsmb'(G)$, we can modify $\cH$ by deleting some winning sets or replacing  an $e\in E(\cH)$ with an $e' \supseteq e$.

\medskip
If $G$ consists of the components $G_1, \dots ,G_\ell$, then $\cH_G$ consists of the components $\cH_{G_1}, \dots \cH_{G_\ell}$. This fact and Proposition~\ref{prop:H-comp} readily imply the following statement.
\begin{proposition}
	\label{prop:comp} 
	If a disconnected graph $G$ consists of the components $G_1, \dots, G_\ell$ and $\gsmb'(G_1) \le \dots  \le \gsmb'(G_\ell)$, then the following statements hold:
	\begin{itemize}
		\item[$(i)$] $\gsmb'(G)=\gsmb'(G_1)$;
		\item[$(ii)$] $\gsmb'(G_1) \le \gsmb(G)\le \gsmb'(G_2)$.	
	\end{itemize}	
\end{proposition}

\begin{proposition} \label{prop:cut}
	Let $v$ be a cut vertex in a connected graph $G$. If $G_1, \dots, G_\ell$ are the components of $G-v$ indexed so that $\gsmb'(G_1)\le  \dots \le \gsmb'(G_\ell)$, then $$\gsmb'(G) \le \gsmb'(G_2)+1.$$
\end{proposition} 
\proof If $\gsmb'(G_2)=\infty$, the inequality is clearly valid. Otherwise, consider the strategy of Staller when she plays $v$ as the first move in the game on $G$. 
By Proposition~\ref{prop:delete-shrink}, the continuation of the game is the same as it would be on $\cH_G \mid v$. As $v$ is a cut vertex in $G$, the hypergraph $\cH_G \mid v$ corresponds to the disjoint union of the closed neighborhood hypergraphs of $G_1, \dots, G_\ell$ supplemented with an extra hyperedge $e=N_G(v)$. Deleting $e$ from $\cH_G \mid v$, we obtain the hypergraph $\cH'$. By  Proposition~\ref{prop:delete}, we have $\w_M^B(\cH_G \mid v) \le \w_M^B(\cH')$. Observe that $\cH'$ contains $\ell$ components, namely $\cH_{G_1}, \dots , \cH_{G_\ell}$.
By Proposition~\ref{prop:comp}, $\w_M^M(\cH')\le \gsmb'(G_2)$ and, therefore, Staller can ensure that she wins after playing at most $1+ \gsmb'(G_2)$ vertices (together with her first move).  
\qed

As a further consequence of Propositions~\ref{prop:delete-shrink} and \ref{prop:delete}, we get the following statement which is analogous to the No-Skip Lemma proved in \cite{gledel-2019} for the MBD games when Dominator's fast winning strategies and the parameters $\gmb(G)$ and $\gmb'(G)$ were studied.

\begin{lemma} [No-Skip Lemma] 
	\label{lem:no-skip} 
	In an optimal strategy of Staller to achieve $\gsmb(G)$ or $\gsmb'(G)$ it is never an advantage for her to skip a move. Moreover, if Dominator skips a move it can never disadvantage Staller. 
\end{lemma}

We close this section with a simple but important relation between the minimum degree $\delta(G)$ and the invariants $\gsmb(G)$, $\gsmb'(G)$.

\begin{proposition}
	\label{prop:mindeg}
	If $G$ is a graph with minimum degree $\delta(G)$, then 
	\begin{equation}
	\label{eq:min-deg}
	\delta(G)+1 \leq \gsmb'(G) \leq \gsmb(G).
	\end{equation}
	\end{proposition}

\proof
Clearly, $\delta(G)+1 \leq \gsmb'(G)$ as Staller can win the game only by claiming a closed neighborhood a vertex.  
To show $\gsmb'(G) \leq \gsmb(G)$, we consider the D-game as an S-game in which Staller skips the first move. By Lemma~\ref{lem:no-skip}, $\gsmb'(G) \leq \gsmb(G)$. \qed

It is clear that the presence of an isolated vertex implies $\gsmb'(G)=1$ and we may easily construct a graph $G_t$ for each $t\ge 1$ such that $\delta(G)=0$, $\gsmb'(G)=1$ and $\gsmb(G)=t$. In Section~\ref{sec:relation}, we will show that for every triple of integers $(r,s,t)$ with $2 \le r \le s \le t$ there exists a graph $G$ satisfying $\delta(G)+1 = r$, $\gsmb(G) = s$, and $\gsmb'(G) = t$.



\section{A realization theorem for $\delta$, $\gsmb'$, and $\gsmb$}
\label{sec:relation}

We start with the definition of two families of graphs. Let $F_1$ be isomorphic to $P_3$ and  $X_1 \subseteq V(F_1)$ contain the two leaves and let $z_1$ denote the support (central) vertex. 
For $k \ge 2$, we construct $F_k$ by using two copies of $F_{k-1}$, namely $F_{k-1}^1$ and $F_{k-1}^2$,  and adding a new vertex $z_k$ which is adjacent to every vertex in $X_{k-1}^1 \cup X_{k-1}^2$. Let us define $X_k=X_{k-1}^1 \cup X_{k-1}^2$ and  $Y_k=V(F_k) \setminus (X_k \cup \{z_k\})$.
 It is clear by the recursive definition  that $|X_k|=2^k$ and each vertex from $X_k$ is of degree $k$. Moreover, $Y_k$ contains $2^{k-i}$ vertices of degree $2^i$ for all $i \in [k-1]$ and, therefore, $|Y_k|=2^k-2$. 

We define  $F'_k$ by modifying $F_k$ as follows. For $i \in [2]$, let $F'_i\cong F_i$.  If $k \ge 3$, $F'_k$ is obtained from $F_k$ by adding a complete graph of order $k-1$. We denote by $Y^+_k$  the set of these $k-1$ new vertices and supplement the graph with a join between $Y^+_k$ and $Y_k$. Observe that $\deg_{F_k'}(v) > k$ holds for all  $v \in Y_k \cup Y_k^+$ and therefore, $\delta(F_k')=k$.

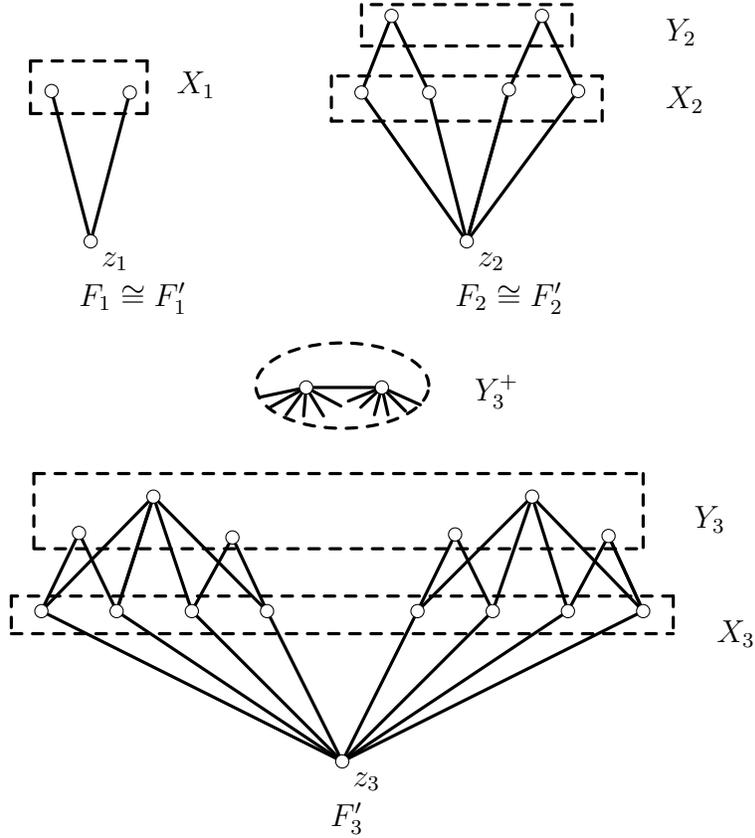
\begin{figure}[t!] 
	\begin{center} 
		
		\definecolor{ffffff}{rgb}{1.,1.,1.}
		\begin{tikzpicture}[line cap=round,line join=round,>=triangle 45,x=1.0cm,y=1.0cm]
		\clip(2.379225210525689,1) rectangle (13.075185726746307,6.990760478820598);
		\draw [line width=1.2pt] (4.,2.)-- (3.48,4.);
		\draw [line width=1.2pt] (4.,2.)-- (4.52,3.98);
		\draw [line width=1.2pt] (9.,2.)-- (8.5,3.98);
		\draw [line width=1.2pt] (9.,2.)-- (7.6,3.98);
		\draw [line width=1.2pt] (9.,2.)-- (9.56,4.02);
		\draw [line width=1.2pt] (9.56,4.02)-- (9.,2.);
		\draw [line width=1.2pt] (9.,2.)-- (10.48,4.);
		\draw [line width=1.2pt] (8.,5.)-- (7.6,3.98);
		\draw [line width=1.2pt] (8.,5.)-- (8.5,3.98);
		\draw [line width=1.2pt] (9.56,4.02)-- (10.,5.);
		\draw [line width=1.2pt] (10.,5.)-- (10.48,4.);
		
		\draw [line width=1.2pt,dash pattern=on 4pt off 4pt] (3.2,4.4)-- (3.2,3.7);
		\draw [line width=1.2pt,dash pattern=on 4pt off 4pt] (3.2,3.7)-- (4.75,3.7);
		\draw [line width=1.2pt,dash pattern=on 4pt off 4pt] (4.75,3.7)-- (4.751,4.4);
		\draw [line width=1.2pt,dash pattern=on 4pt off 4pt] (4.751,4.4)-- (3.2,4.4);
		\draw [line width=1.2pt,dash pattern=on 4pt off 4pt] (7.6,5.15)-- (7.6,4.6);
		\draw [line width=1.2pt,dash pattern=on 4pt off 4pt] (7.6,4.6)-- (10.4,4.6);
		\draw [line width=1.2pt,dash pattern=on 4pt off 4pt] (10.4,4.6)-- (10.4,5.15);
		\draw [line width=1.2pt,dash pattern=on 4pt off 4pt] (10.4,5.15)-- (7.6,5.15);
		\draw [line width=1.2pt,dash pattern=on 4pt off 4pt] (7.2,4.2)-- (7.2,3.6);
		\draw [line width=1.2pt,dash pattern=on 4pt off 4pt] (7.2,3.6)-- (10.8,3.6);
		\draw [line width=1.2pt,dash pattern=on 4pt off 4pt] (10.8,3.6)-- (10.8,4.2);
		\draw [line width=1.2pt,dash pattern=on 4pt off 4pt] (10.8,4.2)-- (7.2,4.2);
		
		\draw (4,2) node[anchor=north west] {$z_1$};
		\draw (3.7,1.6) node[anchor=north west] {$F_1 \cong F'_1$};
		\draw (9,2) node[anchor=north west] {$z_2$};
		\draw (8.7,1.6) node[anchor=north west] {$F_2 \cong F'_2$};
		\draw (11.5,5.1) node[anchor=north west] {$Y_2$};
		\draw (11.5,4.2) node[anchor=north west] {$X_2$};
		\draw (5.8,4.4) node[anchor=north east] {$X_1$};
		
		\begin{scriptsize}
		\draw [fill=ffffff] (4.,2.) circle (2.5pt);
		\draw [fill=ffffff] (3.48,4.) circle (2.5pt);
		\draw [fill=ffffff] (4.52,3.98) circle (2.5pt);
		\draw [fill=ffffff] (9.,2.) circle (2.5pt);
		\draw [fill=ffffff] (8.5,3.98) circle (2.5pt);
		\draw [fill=ffffff] (7.6,3.98) circle (2.5pt);
		\draw [fill=ffffff] (9.56,4.02) circle (2.5pt);
		\draw [fill=ffffff] (10.48,4.) circle (2.5pt);
		\draw [fill=ffffff] (8.,5.) circle (2.5pt);
		\draw [fill=ffffff] (10.,5.) circle (2.5pt);
		
		\end{scriptsize}
		\end{tikzpicture}

		\begin{tikzpicture}[line cap=round,line join=round,>=triangle 45,x=1.0cm,y=1.0cm]
		\clip(3.435616372621552,-6) rectangle (15.336523058107764,0.883499072953899);
		
		\draw [line width=1.2pt] (9.,-5.)-- (5.,-3.);
		\draw [line width=1.2pt] (6.,-3.)-- (9.,-5.);
		\draw [line width=1.2pt] (7.,-3.)-- (9.,-5.);
		\draw [line width=1.2pt] (8.,-3.)-- (9.,-5.);
		\draw [line width=1.2pt] (10.,-3.)-- (9.,-5.);
		\draw [line width=1.2pt] (11.,-3.)-- (9.,-5.);
		\draw [line width=1.2pt] (12.,-3.)-- (9.,-5.);
		\draw [line width=1.2pt] (13.,-3.)-- (9.,-5.);
		\draw [line width=1.2pt] (5.5,-1.96)-- (5.,-3.);
		\draw [line width=1.2pt] (5.5,-1.96)-- (6.,-3.);
		\draw [line width=1.2pt] (7.,-3.)-- (7.54,-2.02);
		\draw [line width=1.2pt] (7.54,-2.02)-- (8.,-3.);
		\draw [line width=1.2pt] (10.5,-1.98)-- (10.,-3.);
		\draw [line width=1.2pt] (10.5,-1.98)-- (11.,-3.);
		\draw [line width=1.2pt] (12.54,-2.)-- (13.,-3.);
		\draw [line width=1.2pt] (13.,-3.)-- (12.54,-2.);
		\draw [line width=1.2pt] (12.54,-2.)-- (12.,-3.);
		\draw [line width=1.2pt] (6.489247075554898,-1.4768749298540385)-- (5.,-3.);
		\draw [line width=1.2pt] (6.489247075554898,-1.4768749298540385)-- (6.,-3.);
		\draw [line width=1.2pt] (6.,-3.)-- (6.489247075554898,-1.4768749298540385);
		\draw [line width=1.2pt] (6.489247075554898,-1.4768749298540385)-- (7.,-3.);
		\draw [line width=1.2pt] (6.489247075554898,-1.4768749298540385)-- (8.,-3.);
		\draw [line width=1.2pt] (11.52361120741799,-1.4768749298540385)-- (11.,-3.);
		\draw [line width=1.2pt] (11.52361120741799,-1.4768749298540385)-- (10.,-3.);
		\draw [line width=1.2pt] (11.52361120741799,-1.4768749298540385)-- (12.,-3.);
		\draw [line width=1.2pt] (11.52361120741799,-1.4768749298540385)-- (13.,-3.);
		\draw [line width=1.2pt] (8.519498840207882,-0.024337081972229144)-- (8.040821594883194,-0.2884348724961945);
		\draw [line width=1.2pt] (8.519498840207882,-0.024337081972229144)-- (8.255401049683917,-0.38747154394268146);
		\draw [line width=1.2pt] (8.519498840207882,-0.024337081972229144)-- (8.486486616392387,-0.4039776558504293);
		\draw [line width=1.2pt] (8.519498840207882,-0.024337081972229144)-- (8.750584406916353,-0.4039776558504293);
		\draw [line width=1.2pt] (8.519498840207882,-0.024337081972229144)-- (7.892266587713465,-0.1563859772342118);
		\draw [line width=1.2pt] (8.519498840207882,-0.024337081972229144)-- (8.981669973624822,-0.2719287605884466);
		\draw [line width=1.2pt] (9.5263716665805,-0.024337081972229144)-- (9.245767764148788,-0.2884348724961945);
		\draw [line width=1.2pt] (9.5263716665805,-0.024337081972229144)-- (9.427334995134014,-0.38747154394268146);
		\draw [line width=1.2pt] (9.5263716665805,-0.024337081972229144)-- (9.60890222611924,-0.38747154394268146);
		\draw [line width=1.2pt] (9.5263716665805,-0.024337081972229144)-- (9.83998779282771,-0.337953208219438);
		\draw [line width=1.2pt] (9.5263716665805,-0.024337081972229144)-- (9.08070664507131,-0.2059043129574553);
		\draw [line width=1.2pt] (9.5263716665805,-0.024337081972229144)-- (10.038061135720683,-0.2554226486806988);
		
		\draw [line width=1.2pt,dash pattern=on 4pt off 4pt] (4.9,-1.17)-- (4.9,-2.17);
		\draw [line width=1.2pt,dash pattern=on 4pt off 4pt] (4.9,-1.17)-- (13.,-1.17);
		\draw [line width=1.2pt,dash pattern=on 4pt off 4pt] (13.,-1.17)-- (13.,-2.17);
		\draw [line width=1.2pt,dash pattern=on 4pt off 4pt] (4.9,-2.17)-- (13.,-2.17);
		\draw [line width=1.2pt,dash pattern=on 4pt off 4pt] (4.6,-2.8)-- (4.6,-3.3);
		\draw [line width=1.2pt,dash pattern=on 4pt off 4pt] (4.6,-3.3)-- (13.4,-3.3);
		\draw [line width=1.2pt,dash pattern=on 4pt off 4pt] (13.4,-2.8)-- (13.4,-3.3);
		\draw [line width=1.2pt,dash pattern=on 4pt off 4pt] (13.4,-2.8)-- (4.6,-2.8);
		\draw [rotate around={0.:(9.,0.)},line width=1.2pt,dash pattern=on 4pt off 4pt] (9.,0.) ellipse (1.1499974196334808cm and 0.5678856092239154cm);
		\draw [line width=1.2pt] (8.519498840207882,-0.024337081972229144)-- (9.5263716665805,-0.024337081972229144);
		
		\draw (13.834466874502708,-2.978931113459094) node[anchor=north west] {$X_3$};
		\draw (13.537356860163246,-1.460368817946293) node[anchor=north west] {$Y_3$};
		\draw (10.615775052491873,0.2892790442749771) node[anchor=north west] {$Y^+_3$};
		\draw (9,-5) node[anchor=north west] {$z_3$};
		\draw (8.7,-5.4) node[anchor=north west] {$F'_3$};
		\begin{scriptsize}
		
		\draw [fill=ffffff] (9.,-5.) circle (2.5pt);
		\draw [fill=ffffff] (5.,-3.) circle (2.5pt);
		\draw [fill=ffffff] (6.,-3.) circle (2.5pt);
		\draw [fill=ffffff] (7.,-3.) circle (2.5pt);
		\draw [fill=ffffff] (8.,-3.) circle (2.5pt);
		\draw [fill=ffffff] (10.,-3.) circle (2.5pt);
		\draw [fill=ffffff] (11.,-3.) circle (2.5pt);
		\draw [fill=ffffff] (12.,-3.) circle (2.5pt);
		\draw [fill=ffffff] (13.,-3.) circle (2.5pt);
		\draw [fill=ffffff] (5.5,-1.96) circle (2.5pt);
		\draw [fill=ffffff] (7.54,-2.02) circle (2.5pt);
		\draw [fill=ffffff] (10.5,-1.98) circle (2.5pt);
		\draw [fill=ffffff] (12.54,-2.) circle (2.5pt);
		\draw [fill=ffffff] (6.489247075554898,-1.4768749298540385) circle (2.5pt);
		\draw [fill=ffffff] (11.52361120741799,-1.4768749298540385) circle (2.5pt);
		\draw [fill=ffffff] (8.519498840207882,-0.024337081972229144) circle (2.5pt);
		\draw [fill=ffffff] (9.5263716665805,-0.024337081972229144) circle (2.5pt);
				\end{scriptsize}
			\end{tikzpicture}
	\end{center}
	\caption{Graphs $F'_1\cong F_1$, $F'_2\cong F_2$, and $F'_3$. In $F_3'$ every vertex in $Y^+_3$ is adjacent to every vertex in $Y_3$.} 
	\label{fig:F'}
\end{figure}

\begin{lemma} \label{F'}
		For every positive integer $k$, $\gsmb'(F'_k)=k+1$ and $\gsmb(F'_k)= \infty$.
\end{lemma}

\proof 
 Let $k$ be a positive integer. Since $\delta(F_k')=k$, Proposition~\ref{prop:mindeg} shows that $\gsmb'(F'_k) \ge k+1$.
 To prove $\gsmb'(F'_k) \le k+1$, we describe an appropriate strategy for Staller in the S-game.
 Staller plays $z_k$ in her first move.
 According to the construction, $F_k'-z_k$ contains two disjoint copies of $F_{k-1}$. Thus, after Dominator's first move, there remains a subgraph $F_{k-1}^i$ such that none of its vertices are played and Staller may claim its center $z_{k-1}^i$ in her second move. Then, $F_{k-1}^i- z_{k-1}^i$ contains two copies of $F_{k-2}$ and at least one of them remains unplayed after the next move of Dominator. From this unplayed $F_{k-2}^j$ subgraph, Staller plays the center $z_{k-2}^j$. Continuing this way, in her $\ell^{\rm th}$ move, $\ell \in [k]$, Staller plays the center of a copy of $F_{k-\ell +1}$. As the described strategy ensures, after Staller's $k^{\rm th}$ move, there remain two vertices, say $x_1$  and $x_2$, such that all vertices from  $N_{F_k'}(x_1)=N_{F_k'}(x_2)$ have been already claimed by Staller. After Dominator's $k^{\rm th}$ move, Staller may play one of these vertices, say $x_1$, and win the game by claiming $N_{F_k'}[x_1]$. This strategy shows that $\gsmb'(F'_k) \le k+1$ and then $\gsmb'(F'_k) = k+1$ follows..
 
 \medskip 
 We next show that $\gsmb(F'_k)= \infty$ for every positive integer $k$. 
 Dominator wins the D-game in $F_1'$ with his first move by playing $z_1$. 
 In $F_2'$, Dominator first plays $z_2$ and then he can claim one (unplayed) vertex from each component of $F_2'-z_2$ as his second and third move. This way Dominator can always claim a dominating set in $F_2'$.  
 If $k\ge 3$,  Dominator may start the game by playing $z_k$. We remark that there remains at least one unplayed vertex in $Y^+_k$ after the first move of Staller. If Dominator responds by playing such a vertex,  he claims a dominating set with this move. We therefore conclude $\gsmb(F'_k)= \infty$. 
\qed

By Proposition~\ref{prop:mindeg}, $\delta(G)+1 \leq \gsmb'(G) \leq \gsmb(G)$ holds for every graph $G$. The following realization theorem implies the sharpness of this inequality (\ref{eq:min-deg}), and also that the differences $\gsmb(G)-\gsmb'(G)$ and $\gsmb'(G)- \delta(G)$ can be arbitrarily large.
  
\begin{theorem}
	\label{thm:mindeg}
	For every three integers $r,s,t$ with $2\le r\le s\le t$, there exists a graph $G$ such that $\delta(G)+1 = r$, $\gsmb'(G) = s$, and $\gsmb(G) = t$.
\end{theorem}

\proof
Let $r, s, t$ be integers such that $2\le r\le s\le t$ and let $G_{r,s,t}$ be the disjoint union of a copy of $K_r$, a copy of $F'_{s-1}$, and a copy of $F'_{t-1}$. We show that $\delta(G_{r,s,t})+1 = r$, $\gsmb'(G_{r,s,t}) = s$, and $\gsmb(G_{r,s,t}) = t$ hold.

Notice first that $\delta(K_r)=r-1$, $\delta(F'_{s-1})=s-1$, and $\delta(F'_{t-1})=t-1$. Thus, $\delta(G_{r,s,t})+1 = r$.
As $r\ge 2$ is supposed, Dominator wins in $K_r$ both in a D-game and S-game that is, $\gsmb(K_r)= \gsmb'(K_r)=\infty$. By Lemma~\ref{F'}, $\gsmb'(F'_{s-1})=s$, $\gsmb'(F'_{t-1})=t$, and $\gsmb(F'_{s-1})=\gsmb(F'_{t-1}) = \infty$.
Proposition~\ref{prop:comp} (i) immediately implies $\gsmb'(G)=s$.

Now consider a D-game played on $G_{r,s,t}$. By Proposition~\ref{prop:comp} (ii), $\gsmb(G_{r,s,t}) \le t$ holds.
 To prove the other direction, we present a strategy for Dominator which guarantees that Staller needs to play at least $t$ vertices to win.
Dominator first plays $z_{s-1}$ in the component $F'_{s-1}$. 
Then, if Staller plays a vertex from a component, in the next turn Dominator replies with an optimal move in the same component. This ensures that Staller cannot win by claiming a closed neighborhood of a vertex from $V(K_r) \cup V(F'_{s-1})$. As $\delta(F'_{t-1})=t-1$, it needs at least $t$ moves to claim a closed neighborhood of a vertex from $F'_{t-1}$. It follows that  $\gsmb(G_{r,s,t}) \ge t$ and we infer $\gsmb(G_{r,s,t}) = t$
\qed

The disconnected graph $G_{r,s,t}$ in the above proof can be replaced with a connected graph that still satisfies the required properties. Indeed, fix a vertex $v \in V(K_r)$ and add the edges $vz_{s-1}$ and $vz_{t-1}$ to $G_{r,s,t}$. As $r \ge 2$, the obtained graph $G'_{r,s,t}$ is of minimum degree $r-1$, and it can be shown that both $\gsmb'(G'_{r,s,t}) = s$ and $\gsmb(G'_{r,s,t}) = t$ remain true.

We note, related to Theorem~\ref{thm:mindeg}, that the statement remains true if the cases $t=\infty$ and $s=t=\infty$  are allowed. Indeed, deleting the component $F'_{t-1}$ from $G_{r,s,t}$, we obtain a graph with $\delta+1=r$, $\gsmb'=s$ and $\gsmb=\infty$, while the deletion of both $F'_{s-1}$ and $F'_{t-1}$ leaves $K_r$ with $\delta(K_r)+1 = r$ and $\gsmb'(K_r)=\gsmb(K_r)=\infty$.

\section{Deleting a leaf and a support vertex} \label{sec:support}
 The deletion of a strong support vertex and all the adjacent leaves from a graph $G$ might increase $\gsmb'(G)$ and $\gsmb(G)$ drastically.  In fact, the difference can be arbitrarily large. In this section we consider the possible changes of the SMBD-numbers when a weak support vertex and its leaf are removed from a graph.
\begin{proposition} \label{prop:G'_S-game}
	Let $G'$ be a graph obtained from $G$ by removing a weak support vertex and the adjacent leaf. Then, the following inequalities hold:
	$$\gsmb'(G)-1 \le \gsmb'(G') \le \gsmb'(G) \qquad {\rm and} \qquad \gsmb(G') \le \gsmb(G).$$
\end{proposition}
\proof
Let $u$ be a weak support vertex in $G$ that is adjacent to the leaf $v$. Then $u$ is a cut vertex in $G$ and $G-u$ is isomorphic to the disjoint union of $G'=G-\{u,v\}$ and the isolated vertex $v$. By Proposition \ref{prop:cut}, $\gsmb'(G)-1 \le\gsmb'(G')$.

 To prove the inequalities $\gsmb'(G') \leq \gsmb'(G)$ and $\gsmb(G') \le \gsmb(G)$, we use the imagination strategy. Let Game 1 be an S-game on  $G'$ where Dominator plays optimally, and let Game 2 be an S-game on $G$ in which Staller applies an optimal strategy. Whenever Dominator plays a vertex in Game 1, he copies this move to Game 2. There Staller replies optimally on $G$, say by playing $s_i$, and copies this move to Game 1 if it is possible. In particular, whenever $s_i \in V(G)\setminus \{u,v\}$, Staller can copy the move to Game 1. Observe that if a winning set in Game 2 contains $v$, it also contains $u$.  We may therefore assume that Staller never plays $v$ while $u$ is also available. If $s_i = u$, then Staller will not play any vertex in Game 1. To continue Game 2, Staller imagines that Dominator replies by playing $v$ in $G$. Then, if Game 2 is not yet over, Staller selects an optimal move $s_{i+1}$ in Game 2 and copies it to Game 1. From this moment, the process continues by simply copying the players' optimal moves to the other game.
 
 For $i \in [2]$, let $t_i$ be the number of Staller's moves until she wins in Game $i$ and set $t_i=\infty$ if she does not win.  
 Observe that  $N_{G'}[x] \subseteq  N_{G}[x]$ holds for every vertex  $x \in V(G)\setminus \{u,v\}$ and, by the described strategy, Staller cannot win in $G$ by claiming the entire $N_{G}[u]$ or $N_{G}[v]$. 
 Thus, when Game 2 ends by the move $s_{t_2}$ and it is copied to Game 1, Staller has already claimed a closed neighborhood in Game 1. 
 We may therefore conclude that $t_1 \le t_2$.
 Suppose first that Staller starts in both Game 1 and Game 2.  Since Dominator plays optimally in Game 1 but Staller might not, we infer $\gsmb'(G') \le t_1$. Similarly, as Dominator might not follow an optimal strategy in Game 2, but Staller plays optimally, we have $ t_2 \le \gsmb'(G)$. This proves $\gsmb'(G') \le \gsmb'(G)$. If Dominator starts in both Game 1 and Game 2, then the same argumentation yields $\gsmb(G') \le \gsmb(G)$.
 \qed
 
 It is not hard to see that $\gsmb(G)-1 \le \gsmb(G')$ is not always true for $G$ and $G'$ under the conditions of Proposition~\ref{prop:G'_S-game}. As an example, we may consider a graph $G$ such that a weak support vertex $u$ is adjacent to all vertices.
 In this graph Dominator may win the D-game with his first move $d_1=u$ that implies $\gsmb(G)= \infty$ even if $\gsmb(G')$ is finite.  
 However we can prove that $\gsmb(G)-1 \le \gsmb(G')$ holds under stronger conditions. Note that the removal of a weak support vertex of degree $2$ and its leaf, had a central role when Dominator's fast winning strategies were studied in \cite{gledel-2019}.
 
 \begin{proposition} \label{prop: G'_D-game}
 Let $G$ be a graph with a weak support vertex of degree $2$ and $G'$ a graph obtained from $G$ by removing this weak support vertex and the adjacent leaf. 
 	Then, $\gsmb(G)-1 \le \gsmb(G')$ holds.
 \end{proposition}
 
 \proof
 Let $u$ be a weak support vertex in $G$ adjacent to the leaf $v$ and to the non-leaf vertex $u'$. We set $G'=G-\{u,v\}$ and note that the closed neighborhoods of vertices in $V(G)\setminus \{u, u', v\}$ are the same in $G$ as in $G'$. We prove  $\gsmb(G)-1 \le \gsmb(G')$ by using the imagination strategy and comparing the winning sets in the two games after a vertex from $\{u,v,u'\}$ was played. Let Game 1 be a D-game on $G$ where Dominator plays optimally, and Game 2 be a D-game on $G'$ where Staller plays optimally. If neither of $u$, $u'$ and $v$ has been played until a moment in the game, then the moves of the players are simply copied to the another game. That is, whenever Dominator plays a vertex $d_i$ in Game 1, he copies this move to Game 2. There Staller replies optimally on $G'$, say by playing $s_i$, and copies this move to Game 1. Under the given conditions, it is always possible. If a vertex from $\{u,v,u'\}$ is played, the corresponding move(s) in the other game will be specified in Case 1 and 2.
  For $i \in [2]$, let $t_i$ be the number of Staller's moves until she wins in Game $i$ and set $t_i=\infty$ if she does not win. We prove $t_1 \le t_2+1$ by considering three cases.
    
  \paragraph{Case 1.}
  If Dominator plays first a vertex $d_i$ from $\{u,v,u'\}$, this move belongs to Game 1. We may suppose that Dominator's corresponding move in Game 2 is $u'$. After these moves, each remaining winning set of Staller in Game 2 is also a winning set in Game 1. By Proposition~\ref{prop:delete} (i), this proves $t_1 \le t_2$.
    
  \paragraph{Case 2.} 
  If Staller plays first a vertex $s_i$ from $\{u,v,u'\}$, then $s_i=u'$ as this move belongs to Game 2.
 Note that Case 1 excludes Case 2 in a game and, therefore, $u$, $u'$ and $v$ all are unplayed in Game 1.  We may set  $s_i= u$ in Game 1. Dominator, who plays optimally there, has to reply by claiming $v$ as otherwise Staller could win in her next move by playing $v$. These two moves $u$ and $v$ are not copied to Game 2. Instead, we set $s_{i+1}=u'$ in Game 1 and observe that the winning sets in Game 1 and Game 2 are the same, while the number of moves of Staller in Game 1 is one more than the number of her moves in Game 2. We may infer $t_1=t_2+1$. Remark that this argumentation is also true if Staller finishes Game 2 with the move $s_{t_2}=u'$. 

 \paragraph{Case 3.}
 If neither of $u$, $u'$ and $v$ is played during the games, then the moves are just copied from one game to the other. Further, under this condition, Staller cannot win by claiming the entire $N_{G'}[u']$, $N_{G}[u']$, $N_{G}[u]$, or $N_{G}[v]$. For the remaining vertices, we have $N_{G}[x]=N_{G'}[x]$ and thus, the two games end together. This shows $t_1=t_2$.

 \vspace*{0.5cm}

 We may therefore conclude that $t_1 \le t_2+1$. Our argumentation is also valid for the cases when Staller does not win. For this situation, we suppose that every vertex is played in $G$ and $G'$ but no winning set is claimed by Staller.
 
   Since Dominator plays optimally all along Game~1, we infer $\gsmb(G) \le t_1$. Similarly, as Staller follows an optimal strategy in Game~2, we have $ t_2 \le \gsmb(G')$. Together with $t_1 \le t_2+1$ we conclude $\gsmb(G) \le \gsmb(G')+1$. 
\qed

\section{Paths}
\label{sec:path}
The ``pairing strategy'' of Maker \cite{hefetz}  ensures that he can win in the Maker-Breaker game if an approriate set of vertex pairs can be defined. Its immediate consequence for the MBD game shows that Dominator has a winning strategy in both the D-game and S-game if the graph admits a perfect matching. Here we state the lemma in a more general form.
\begin{lemma}
	\label{lem:pairing}
	Consider an MBD game on $G$ and let $X$ and $Y$ be the set of vertices played by Dominator and Staller, respectively, until a moment during the game. If there exists a matching $M$ in $G-(X\cup Y)$ such that $V(G) \setminus V(M) \subseteq N_G[X]$,
	then Dominator has a strategy to win the continuation of the game, no matter who plays the next vertex.
\end{lemma} 
\proof By the conditions of the lemma,
$$(X \cup Y)\subseteq (V(G) \setminus V(M)) \subseteq N_G[X]$$ 
that is all vertices from $X\cup Y$ are dominated by the first $|X|$ moves of Dominator.
In the continuation Dominator can win the game by applying the following strategy. If Staller plays a vertex $v$ such that $uv \in E(M)$ and $u$ has not been played yet, then Dominator plays $u$ in his next turn. Otherwise, Dominator plays an arbitrary unplayed vertex. As $V(M)$ does not contain any vertices from $Y$, this startegy  makes sure that Dominator plays at least one vertex from each edge of the matching (unless the game finishes earlier as the entire $V(M)$ is dominated by Dominator's moves). We conclude that Dominator's moves dominate every vertex from $V(G)= V(M) \cup (V(G)\setminus V(M))$ and hence, he wins the game.  
\qed
\medskip

As it was shown in \cite{duchene-2019+}, Dominator can win the MBD game on all paths and cycles if he is the first player. Moreover, Dominator also wins the S-game if it is played on a cycle or a path of an even order. On the other hand, if $n$ is odd, $\gsmb'(P_n) < \infty$. 

By Proposition~\ref{prop:G'_S-game},   $\gsmb'(P_{n+2}) \le \gsmb'(P_n)+1$. Starting with $\gsmb'(P_1)=1$ and $\gsmb'(P_3)=2$, this yields  $\gsmb'(P_{2k+1}) \le  k+1$. For paths of even order the corresponding invariant satisfies $\gmb'(P_{2k})=k$ (see \cite{gledel-2019}). These facts might suggest that $\gsmb'(P_{2k+1})$ is equal to or not far from $k+1$. However, as we will show, the value of  $\gsmb'(P_{2k+1})$ is significantly different from $k+1$.
The exact formula is established in the following theorem. If $n \ge 3$ and $n$ is odd, we may also cite the statement of the theorem in the form of  $\gsmb'(P_n)=\lceil \log_2 n \rceil$.

\begin{theorem} \label{thm:path}
 If $n$ is an odd positive integer, then
 \begin{align} \label{form:path}
 \gsmb'(P_n)=\left\lfloor \log_2n \right\rfloor +1.
 \end{align}
 Moreover, Staller has an optimal strategy in the S-game on $P_n$ such that she wins on the closed neighborhood of a vertex $v$ that is at an even distance from the ends of the path.
\end{theorem}
\proof First we observe that $\gsmb'(P_1)=1$ and $\gsmb'(P_3)=2$ and thus, the formula is valid for $n=1$ and $n=3$. Then we proceed by induction on $n$. Suppose that $n \ge 5$ and (\ref{form:path}) is true for every odd integer which is smaller than $n$. We first prove $\gsmb'(P_n) \le  \left\lfloor \log_2n \right\rfloor +1$ by considering two cases.

\paragraph{Upper bound, Case 1.} Assume that $n=4k+3$ for an integer $k \ge 1$ and consider the path $P_n: v_1\dots v_{4k+3}$. The vertex  $v_{2k+2}$ is a cut vertex and both components of $G-v_{2k+2}$ are isomorphic to $P_{2k+1}$ where $2k+1=\frac{n-1}{2}$. Proposition~\ref{prop:cut} and our hypothesis imply
$$\gsmb'(P_n) \le \gsmb'\left(P_{\frac{n-1}{2}}\right)+1= \left\lfloor \log_2\frac{n-1}{2} \right\rfloor +2 = \left\lfloor \log_2 (n-1) \right\rfloor +1. 
$$
To finish the proof we add that, since $n=4k+3$ is not a power of $2$, the equality $\left\lfloor \log_2 (n-1) \right\rfloor = \left\lfloor \log_2 n \right\rfloor$ necessarily holds.

\paragraph{Upper bound, Case 2.} Assume that $n=4k+1$ for an integer $k \ge 1$ and consider the path $P_n: v_1\dots v_{4k+1}$. Now, we choose the cut vertex $v_{2k}$ and consider the components $P^1:v_1\dots v_{2k-1}$ and $P^2: v_{2k+1}\dots v_{4k+1}$ of $P_n -v_{2k}$. Remark that $P^2$ is a path of order $\frac{n+1}{2}$. As  $\gsmb'(P^1) \le \gsmb'(P^2)$, Proposition~\ref{prop:cut} and our hypothesis together imply
$$\gsmb'(P_n) \le \gsmb'\left(P_{\frac{n+1}{2}}\right)+1= \left\lfloor \log_2\frac{n+1}{2} \right\rfloor +2 = \left\lfloor \log_2 (n+1) \right\rfloor +1. 
$$
As $n+1=4k+2$ and $n+1 \ge 6$, it cannot be a power of $2$ and therefore, we get $\left\lfloor \log_2 (n+1) \right\rfloor =\left\lfloor \log_2 n \right\rfloor$. It finishes the proof for the upper bound.

\paragraph{Lower bound.} We now prove that Dominator can ensure that Staller does not win in her first $\left\lfloor \log_2 n \right\rfloor$ moves on the path $P_n: v_1 \dots v_n$. We analyze different cases according to Staller's first move $s_1$.
\begin{itemize}
\item Assume that $s_1=v_i$ for an odd integer $i\in \{1,3,\dots ,n\}$. By symmetry and the condition $n \ge 5$, we may suppose $i \ge 3$. Dominator then replies by playing $d_1=v_{i-1}$. The remaining vertices, except $v_{i-2}$, are covered by the following matching in $G-\{v_i, v_{i-1}\}$: 
$$ v_1v_2, \dots , v_{i-4}v_{i-3}, \, v_{i+1}v_{i+2}, \dots , v_{n-1}v_{n}.
$$
As it satisfies the conditions in Lemma~\ref{lem:pairing}, we may conclude, by Lemma~\ref{lem:pairing}, that Dominator can win the game if Staller's first move is $v_i$ with an odd index $i$. 
\item Suppose now that $s_1=v_i$ and $i$ is even. We may also assume, without loss of generality, that $i\le \frac{n+1}{2}$. Dominator's response is the move $d_1=v_{i-1}$. After these two moves, by Proposition~\ref{prop:delete-shrink}, the winning sets for Staller are exactly the edges of the hypergraph $\cH'=(\cH_{P_n}\! \mid \! v_i)-v_{i-1}$. If $i=2$, then $\cH'$ is isomorphic to the closed neighborhood hypergraph of the path $P_{n-2}$ and, by the induction hypothesis, Staller needs at least $\lfloor \log_2(n-2)\rfloor +1$ further moves there to win. Note that, as $n \ge 5$, the inequality $\lfloor \log_2(n-2)\rfloor +2 \ge \lfloor \log_2n \rfloor +1$ is valid and completes the proof for the case of $i=2$.
 If $i\ge 4$, then $\cH'$ contains two components. One of them is exactly the closed neighborhood hypergraph of the odd path  $P^2: v_{i+1} \dots v_n$. Adding the extra (hyper)edge $\{v_{i-2},v_{i-3}\}$ to the other component of $\cH'$, we obtain the closed neighborhood hypergraph of the even path $P^1: v_1\dots v_{i-2}$. Since $\gsmb'(P^1)=\infty$, Proposition~\ref{prop:comp} and the induction hypothesis together imply that Staller needs at least $\gsmb'(P^2)= \lfloor \log_2(n-i)\rfloor +1$ moves to win the S-game on $P^1\cup P^2$. Since  $\cH'$ can be obtained from the closed neighborhood hypergraph  $\cH_{P^1\cup P^2}$ by deleting the edge $\{v_{i-2},v_{i-3}\}$, Staller needs at least as many moves to win on $\cH'$ as on $\cH_{P^1\cup P^2}$. 
 We may conclude that, if Dominator plays optimally and Staller wins with her  $t^{\rm th}$ move, then
$$ t \ge \left\lfloor \log_2(n-i)\right\rfloor +2 \ge \left\lfloor \log_2\frac{n-1}{2}\right\rfloor +2 =\left\lfloor \log_2(n-1)\right\rfloor +1=\left\lfloor \log_2n\right\rfloor +1,$$
where the last equality is true because $n$ cannot be a power of $2$.  
\end{itemize}
This establishes the lower bound $\gsmb'(P_n)\ge \left\lfloor \log_2n\right\rfloor +1$ and, together with the first part of the proof, we obtain $\gsmb'(P_n)= \left\lfloor \log_2n\right\rfloor +1$ for every odd integer $n$. Analizing Staller's strategy described in the proof of the upper bound, we also infer that the claimed winning set is always one of $N[v_1], N[v_3], \dots, N[v_n]$. 
\qed

\section{Tadpole graphs}
\label{sec:tadpole}

The \emph{tadpole graph} $T(n,k)$ with parameters $n \ge 3$ and $k \ge 1$  is obtained from an $n$-cycle and a (vertex disjoint) path of order $k$ such that one vertex from the cycle and one end of the path are made adjacent. These graphs have been recently studied in \cite{bujtas-2021+} concerning Rall's $1/2$-conjecture \cite{james-2019} on the domination game.

First we prove a short technical lemma.
\begin{lemma} \label{claim:2}
	For every two integers $a$ and $b$ with $a \ge 1$ and $b \ge 2$, the following is true:
	$$\max\{\lfloor \log_2 a\rfloor +1, \lceil \log_2(b-1) \rceil\} \ge  \lceil \log_2(a+b) \rceil -1.
	$$
\end{lemma}
\textit{Proof. } Let $\lceil \log_2(a+b) \rceil=s$ and note that $s\ge 2$.  Suppose for a contradiction that $\lfloor \log_2 a\rfloor +1 \le s-2$ and $ \lceil \log_2(b-1) \rceil \le s-2$ simultaneously hold. These assumptions imply $a < 2^{s-2}$ and $b \le 2^{s-2} +1$. We infer $a+b < 2^{s-1}+1$ and, as $a+b$ is an integer, we get $a+b \le 2^{s-1}$ that contradicts the definition of $s$. \qed
\bigskip

To state the main theorem of this section, we introduce the notation
$$\sigma(n)= \frac{2^{\lceil \log_2 n \rceil}-n}{2}+1.$$

\begin{theorem}
	\label{thm:tadpole}
For every tadpole graph $T(n,k)$ it is true that $\gsmb(T(n,k)) = \infty$ and
$$\gsmb'(T(n,k)) = \begin{cases} 
\lceil \log_2(n+k+\sigma(n))\rceil &  \text{ if $n$ is even and $k$ is odd};\\ 
\infty & \text{otherwise}. 
\end{cases}$$ 
\end{theorem}
\proof Let $v_0v_1\dots v_{n-1}v_0$  and $u_1\dots u_k$ be the cycle and the path of the tadpole graph such that $v_0$ and $u_k$ are adjacent. 

We first prove that Dominator can win the D-game that is $\gsmb(T(n,k)) = \infty$ holds for every pair $(n,k)$ of positive integers with $n \ge 3$. Let Dominator play $d_1=v_0$ as his first move. The edges
$$v_1v_2,\dots , v_{2\lceil n/2 \rceil -3} v_{2\lceil n/2 \rceil -2}, \, u_1u_2, \dots, u_{2\lfloor k/2 \rfloor -1} u_{2\lfloor k/2 \rfloor} $$
form a matching in $T(n,k)-v_0$ that covers all vertices which are not dominated by $v_0$.   Lemma~\ref{lem:pairing} thus implies  $\gsmb(T(n,k)) = \infty$.

Next, we consider the S-game. If both $n$ and $k$ are odd or both are even, then $T(n,k)$ has a perfect matching and Dominator can win by applying the pairing strategy. If $n$ is odd and $k$ is even, we may iteratively remove a weak support vertex and its leaf until a cycle $C_n$ remains. Thus, by Proposition~\ref{prop:G'_S-game}, $\gsmb'(C_n) \le \gsmb'(T(n,k))$ is valid and then, $\gsmb'(C_n)=\infty$ \cite{duchene-2019+} implies $\gsmb'(T(n,k))=\infty$.

\medskip

In the rest of the proof, we suppose that $n$ is even, $k$ is odd, and an S-game is played on $T(n,k)$. From now on, $\ell$ denotes $\lceil \log_2 n \rceil$. Note that $\ell \ge 3$ unless $n=4$.

\paragraph{Upper bound} We consider three cases and show that Staller can win by playing at most $\lceil \log_2(n+k+\sigma(n))\rceil$ vertices. The first two cases together settle the proof for every tadpole graph $T(n,1)$, $n \ge 4$. Then, in Case 3, we proceed by induction on $k$ for every fixed $n$.

\paragraph{Case 1.} $n+k+\sigma(n) \le 2^\ell$.\\
To prove $\gsmb'(T(n,k)) \le \ell$, consider the following strategy of Staller. First, Staller plays $s_1=v_{i}$ with $i=2\lceil \frac{n}{4} \rceil$. Observe that $v_i$ and $v_0$ split the cycle into two nearly equal odd paths that are $P^1: v_1\dots v_{i-1}$ and $P^2: v_{i+1}\dots v_{n-1}$. If $n=4p$, both paths are of order $\frac{n}{2}-1$, and if $n=4p+2$, then $P^1$ is of order $\frac{n}{2}$ and the order of $P^2$ is $\frac{n}{2}-2$. Note that $2 \le i \le n-2$.

In the continuation, if Dominator's first move $d_1$ is not a neighbor of $v_i$, Staller selects her next move  $s_2$ from $\{v_{i-1}, v_{i+1}\}$ such that the distance between $s_2$ and $d_1$ is at least $3$. It can be done as our present condition $n+k+\sigma(n) \le 2^\ell$ excludes $n=4$. Without loss of generality, we may suppose that $s_2=v_{i-1}$. Then, after Dominator's move $d_2$, one of $v_{i-2}$ and $v_{i+1}$ remains unplayed. Playing this vertex as $s_3$, Staller can win the game by claiming $N[v_i]$ or $N[v_{i-1}]$. In this case, Staller wins the game in $3 \le \ell$ moves.

If Dominator plays a neighbor of $s_1$ as his first move $d_1$, the game continues with the winning sets in $(\cH_{T(n,k)}\mid s_1) -d_1$. Suppose that $d_1=v_{i+1}$. Deleting the winning sets $N[v_{i+3}], \dots, N[v_{n-1}]$ from $(\cH_{T(n,k)}\mid s_1) -d_1$, we obtain the hypergraph $\cH'$. By Proposition~\ref{prop:delete}, $\w_M^M(\cH') \ge \w_M^M((\cH_{T(n,k)}\mid s_1) -d_1)$. Note that $\cH'$ can be obtained from the closed neighborhood hypergraph of the path $P': v_{i-1}\dots v_0u_k\dots u_1$ by replacing the edge $\{v_1,v_0,u_k\}$ with $\{v_1,v_0,u_k, v_{n-1}\}$.  Then, Staller continues playing according to the optimal strategy on the path $P'$ that makes sure that the claimed winning set is either $v_j$ or $u_j$ with an odd index $j$. By Theorem~\ref{thm:path}, such a winning strategy exists and it can be applied for $\cH'$.  The order of $P'$ is at most $\frac{n}{2}+1+k$, and therefore, Staller can win the game in at most $1+ \lceil \log_2(\frac{n}{2}+1+k)\rceil$ moves. We observe that
$$1+ \left\lceil \log_2\left( \frac{n}{2}+1+k \right) \right\rceil = 1+ \lceil \log_2(n+k+ \sigma(n)-2^{\ell-1}) \rceil \le 1+\lceil \log_2 (2^\ell- 2^{\ell-1}) \rceil= \ell,
$$
where the inequality follows from the condition $n+k+\sigma(n) \le 2^\ell$. 

 If $d_1=v_{i-1}$, the argumentation is similar. Here we consider the path $P': v_{i+1}\dots v_{n-1}v_0u_k\dots u_1$ which is of order $n-i +k \le \frac{n}{2}+k$.

We conclude that, under the present condition, Staller has a strategy to win the S-game by playing at most $\ell$ vertices. Therefore, $\gsmb'(T(n,k)) \le \lceil \log_2(n+k+\sigma(n))\rceil$ holds if the parameters $n$ and $k$ satisfy the condition in Case~1.

\paragraph{Case 2.} $n \ge 2^\ell-2$   and $k=1$.\\
As $\sigma(n)= 1$ if $n =2^\ell$, and $\sigma(n)=2$ if $n =2^\ell-2$, we have $\lceil \log_2(n+k+\sigma(n))\rceil = \ell+1$ under the present condition.	
To show $\gsmb'(T(n,1)) \le \ell+1$, we apply Proposition~\ref{prop:cut} for the cut vertex $v_0$. The components of $T(n,1)-v_0$ are an isolated vertex with $\gsmb'(P_1)=1$ and an $(n-1)$-path with $\gsmb'(P_{n-1})=\lceil \log_2(n-1) \rceil=\ell.$  Proposition~\ref{prop:cut} therefore implies $\gsmb'(T(n,1)) \le \ell+1$.

\paragraph{Case 3.} $2^{s-1} < n+k+\sigma(n) \le 2^{s}$ and $s \ge \ell+1$.\\
To prove $\gsmb'(T(n,k)) \le s$ we identify appropriate cut vertices again. If $k \le 2^{s-1}$, we consider $T(n,k)-v_0$. The two components are isomorphic to two paths, one is of order $n-1$ and the other is of order $k$. by definition of $\ell$, we have $n-1 < 2^\ell \le 2^{s-1}$ that gives  $\max\{n-1,k\} \le 2^{s-1}$ in turn and we conclude, by Proposition~\ref{prop:cut} that $\gsmb'(T(n,k))\le (s-1)+1=s$. 

If $k >2^{s-1}$, we choose the cut vertex $u_{2^{s-1}}$ and consider the two components obtained after the removal of $u_{2^{s-1}}$. One component is a path with $\gsmb'(P_{2^{s-1}-1})= s-1$ while the other component is a tadpole graph $T'=T(n, k-2^{s-1})$. Since  $n+k+\sigma(n) \le 2^{s}$ implies  $n+k- 2^{s-1}+\sigma(n) \le 2^{s-1}$, the induction hypothesis yields $\gsmb'(T')\le s-1$. By Proposition~\ref{prop:cut}, we infer $\gsmb'(T(n,k))\le (s-1)+1=s$ again.

\paragraph{Lower bound.} In the second part, we show that 
$$\gsmb'(T(n,k))\ge \lceil \log_2(n+k+\sigma(n))\rceil$$
 holds for every pair $(n,k)$ with $n\ge 4$ even and $k \ge 1$ odd. First we eliminate some worst first moves of Staller which make possible Dominator's winning in the continuation. Let $s_1$ and $d_1$ denote Staller's and Dominator's first move in the S-game on $T(n,k)$.
 
\paragraph{Case 1}
 $s_1=u_i$ and $i$ is odd or $s_1=v_j$ and $j$ is odd.\\
 These are Staller's worst first moves. We prove that, under this condition, Dominator has a winning strategy in the continuation of the game.
 \begin{itemize}
 	\item If $s_1=u_1$ and $k=1$, Dominator replies by playing $d_1=v_0$. Then, the graph $G-\{u_1,v_0,v_1\}$ admits a perfect matching that satisfies the conditions of Lemma~\ref{lem:pairing} with $X=\{d_1\}$, $Y=\{s_1\}$ and hence, Dominator can win the continuation of the game.  
 	\item If $s_1=u_1$ and $k \ge 3$, Dominator plays $d_1=u_2$. Here, we consider the perfect matching in $G-\{u_1,u_2,u_3\}$ and conclude, by Lemma~\ref{lem:pairing}, that Dominator can win.
 	\item If $s_1=u_i$ and $i$ is an odd integer with $i \ge 3$, Dominator replies with $d_1=u_{i-1}$. The graph $G-\{u_i,u_{i-1},u_{i-2}\}$ admits a perfect matching and, by applying Lemma~\ref{lem:pairing}, we infer that Dominator can win in the continuation of the game.
 	\item If $s_1=v_1$, Dominator plays $d_1=v_2$. As $G-\{v_1,v_2,v_3\}$ admits a perfect matching, we conclude, by Lemma~\ref{lem:pairing}, that Dominator can win.
 	\item If $s_1=v_j$ and $j$ is an odd integer with $j \ge 3$, Dominator plays $d_1=v_{j-1}$. Since the conditions of Lemma~\ref{lem:pairing} are satisfied by the perfect matching in    $G-\{v_j,v_{j-1},v_{j-2}\}$, Dominator can win the game.  \end{itemize}
 
 \paragraph{Case 2.} $s_1=v_0$.\\
 If $k \le n-1$, Dominator replies by playing $d_1=u_k$ and the continuation of the game reduces to a Maker-Breaker game on $\cH'=(\cH_{T(n,k)}\mid v_0)-u_k$. If $k=1$, then $\cH'$ is the closed neighborhood hypergraph of the path $v_1,\dots, v_{n-1}$. If $k\ge 3$ then, by adding $\{u_{k-1}, u_{k-2}\}$ to the winning sets, we obtain a hypergraph $\cH''$ that consists of two components. $\cH''$  is the closed neighborhood hypergraph of a disjoint union of a path of order $n-1$ and a path of order $k-1$. As $k-1$ is even,  $\gsmb'(P_{k-1})=\infty$. Applying Proposition~\ref{prop:delete}, \ref{prop:comp}, and Theorem~\ref{thm:path} we get
 $$\w_M^M(\cH')\ge \w_M^M(\cH'')=\gsmb'(P_{n-1})= \lceil \log(n-1)\rceil =\ell, 
 $$
 no matter whether $k=1$ or $k \ge 3$.
  Dominator thus can force Staller to play at least $\ell+1$ vertices in the game. On the other hand, the condition $k \le n-1 \le 2^\ell-1$ implies 
  $$ n+k+\sigma(n)= \frac{n}{2}+k+2^{\ell-1}+1 \le 2^{\ell-1}+2^\ell -1 +2^{\ell-1}+1=2^{\ell+1}.
  $$
  Consequently, Dominator can ensure that Staller plays at least $\lceil \log_2(n+k+\sigma(n))\rceil$ vertices.
  
  Assume now that $k \ge n$. Then, Dominator replies with $d_1=v_1$ and the game continues on $\cH'=(\cH_{T(n,k)}\mid v_0)-v_1$. Adding the winning set $\{v_2,v_3\}$ to $\cH'$, we obtain $\cH''$ that is the closed neighborhood hypergraph of two disjoint paths one of them is an even path of order $n-2$ and the other one is of order $k$. By the same reasoning as before, we get that Staller needs to play at least $1+ \lceil \log_2(k) \rceil$ vertices to win. If $\ell+1 \ge \lceil \log_2(n+k+\sigma(n))\rceil$, this already finishes the proof as $\lceil \log_2(k) \rceil \ge  \lceil \log_2(n) \rceil =\ell$. In the other case, $\lceil \log_2(n+k+\sigma(n))\rceil=s \ge \ell+2$ and 
  $$ 2^{s-1} < n+k+\sigma(n)=\frac{n}{2}+k+2^{\ell-1}+1. $$ 
  This, together with the condition $n/2 \le 2^{\ell-1}$ implies $\lceil \log_2(k) \rceil \ge s-1$ and therefore, $1+ \lceil \log_2(k) \rceil \ge s$ as required.

  \paragraph{Case 3.} $s_1=v_i$, $i$ is even and $i\ge 2$.\\  
  By symmetry, we may suppose that $i \ge n/2$ if $n \equiv 0 \pmod{4}$ and  $i \ge n/2+1$ if $n \equiv 2 \pmod{4}$.  Let Dominator reply by playing $d_1=v_{i+1}$. By Proposition~\ref{prop:delete-shrink}, the game continues on $\cH=(\cH_{T(n,k)}\mid v_i) -v_{i+1}$. We modify $\cH$ by adding the winning set $\{v_{i+2},v_{i+3}\}$ (if $i+1\neq n-1$) and replacing $N[v_0]$ and $N[v_{n-1}]$ with $N[v_0]\setminus \{v_{n-1}\}$ and $N[v_{n-1}]\setminus \{v_{0}\}$. The obtained hypergraph is denoted by $\cH'$.  By Proposition~\ref{prop:delete}, $\w_M^M(\cH') \le w_M^M(\cH)$. Observe that $\cH'$ is the closed neighborhood hypergraph of the graph $G'$ that consists of two (or only one) path components that are $P^1: v_{i-1}\dots v_0 u_k  \dots u_1$ and $P^2: v_{i+2}\dots v_{n-1}$. The component $P^1$ is a path of order $i+k$ and  $P^2$ is a path of even order. By Proposition~\ref{prop:comp}, $\gsmb'(G')=\gsmb'(P^1)=\lceil \log_2(i+k) \rceil$. Hence, if $t$ is the number of moves Staller needs to win when $s_1=v_i$ is fixed, then 
  \begin{equation} \label{eq:4}
  t \ge 1+ \w_M^M(\cH') = 1+ \gsmb'(P^1)=1+ \lceil \log_2(i+k) \rceil. 	
  \end{equation}
  Let $s=\lceil \log_2(n+k+\sigma(n))\rceil$. 
  \begin{itemize}
  \item  If $s=\ell$, then (\ref{eq:4}) together with $i+k \ge \frac{n}{2}+1 > 2^{\ell-2}$ gives $t \ge s$.
  \item  If $s=\ell+1$, we consider the inequality
    \begin{equation} \label{eq:5}
  \frac{n}{2}+k =n+k+\sigma(n)-2^{\ell-1}-1 >2^\ell-2^{\ell-1}-1= 2^{\ell-1}-1.
  \end{equation}
  If $n \equiv 0 \pmod{4}$, then $i \ge n/2$ and, as $n/2 +k$ is odd, the strict inequality in (\ref{eq:5}) reduces to $\frac{n}{2}+k \ge 2^{\ell-1}+1$. Then, (\ref{eq:4}) and (\ref{eq:5}) together give $t\ge \ell+1$.
  If $n \equiv 2 \pmod{4}$, then $i \ge 1+ n/2$. By (\ref{eq:5}), we have
  $$i+k \ge \frac{n}{2}+k+1 > 2^{\ell-1}$$
  and, by using (\ref{eq:4}), we conclude $t\ge s$ again.
  \item If $s \ge \ell+2$, then $i+k > 2^{s-2}$. Then, $1+ \lceil \log_2(i+k) \rceil \ge s$ and (\ref{eq:4}) imply the desired inequality $t \ge s$.  
  \end{itemize}
Remark that Cases 1--3 cover all possibilities if $k=1$. Hence, we may proceed by induction on $k$ when the last case is considered. 
  
\paragraph{Case 4.} $s_1=u_i$ and $i$ is even.\\
Applying  Lemma~\ref{claim:2} with $a=i-1$ and $b=n+k-i+\sigma(n)+1$ we obtain 
\begin{equation} \label{eq:6}
\begin{split}
m & =  \max\{\lfloor \log_2(i-1) \rfloor +1, \lceil \log_2(n+k-i+\sigma(n)) \rceil\}\\
& \ge 
\lceil \log_2(n+k+\sigma(n)) \rceil -1.
\end{split}
\end{equation}
If $m=\lfloor \log_2(i-1) \rfloor +1$, Dominator's reply to $s_1=u_i$ is $d_1=u_{i+1}$. Then, the game continues on $\cH'=(\cH_{T(n,k)}\mid u_i)-u_{i+1}$. Adding a new winning set $\{u_{i+2}, u_{i+3}\}$ to $\cH'$, we obtain $\cH''$. Note that playing on $\cH''$ is at least as advantageous for Staller as on $\cH'$. One component of $\cH''$ is the closed neighborhood hypergraph of a tadpole graph $T'=T(n, k-i-1)$. As $k-i-1$ is even, $\gsmb'(T')=\infty$. The other component of $\cH''$ is the closed neighborhood hypergraph of a path of order $i-1$. By Proposition~\ref{prop:comp}, Staller has to play at least $m$ further vertices to win. As (\ref{eq:6}) shows, the total number of Staller's moves is at least $m+1 \ge \lceil \log_2(n+k+\sigma(n)) \rceil$.

The proof is similar if $m=\lceil \log_2(n+k-i+\sigma(n)) \rceil$. Then, Dominator plays $d_1=u_{i-1}$ and then, assuming optimal strategies in the continuation, the game is played on the closed neighborhood hypergraph of a tadpole graph $T(n, k-i) $ and Staller has to play at least $m+1= \lceil \log_2(n+k-i+\sigma(n)) \rceil +1$ vertices to win. By (\ref{eq:6}), it is at least $\lceil \log_2(n+k+\sigma(n)) \rceil$. \qed

\section{Concluding remarks} \label{sec:concluding}

\begin{figure}[t] 
	\begin{center}
		\definecolor{ffffff}{rgb}{1.,1.,1.}
		\begin{tikzpicture}[line cap=round,line join=round,>=triangle 45,x=1.0cm,y=1.0cm]
		\clip(0.7,-0.56) rectangle (15.36,4.1);
		\draw [line width=1.2pt] (2.,3.)-- (2.,2.);
		\draw [line width=1.2pt] (3.,3.)-- (3.,2.);
		\draw [line width=1.2pt] (6.,3.)-- (6.,2.);
		\draw [line width=1.2pt] (4.,1.)-- (2.,2.);
		\draw [line width=1.2pt] (4.,1.)-- (3.,2.);
		\draw [line width=1.2pt] (3.,2.)-- (4.,1.);
		\draw [line width=1.2pt] (4.,1.)-- (6.,2.);
		\draw [line width=1.2pt] (5.,3.)-- (5.,2.);
		\draw [line width=1.2pt] (4.,1.)-- (5.,2.);
		\draw [line width=1.2pt] (10.,3.)-- (10.,2.);
		\draw [line width=1.2pt] (11.,3.)-- (11.,2.);
		\draw [line width=1.2pt] (13.,3.)-- (13.,2.);
		\draw [line width=1.2pt] (14.,3.)-- (14.,2.);
		\draw [line width=1.2pt] (11.54,1.)-- (10.,2.);
		\draw [line width=1.2pt] (11.54,1.)-- (11.,2.);
		\draw [line width=1.2pt] (11.,2.)-- (11.54,1.);
		\draw [line width=1.2pt] (13.,2.)-- (11.54,1.);
		\draw [line width=1.2pt] (14.,2.)-- (11.54,1.);
		\draw [line width=1.2pt] (12.54,1.)-- (14.,2.);
		\draw [line width=1.2pt] (13.,2.)-- (12.54,1.);
		\draw [line width=1.2pt] (12.54,1.)-- (11.,2.);
		\draw [line width=1.2pt] (12.54,1.)-- (10.,2.);
		\draw (4.,0.82) node[anchor=north west] {$S_4^1$};
		\draw (12.16,0.82) node[anchor=north west] {$S_4^2$};
		\begin{scriptsize}
		\draw [fill=ffffff] (2.,3.) circle (2.5pt);
		\draw [fill=ffffff] (2.,2.) circle (2.5pt);
		\draw [fill=ffffff] (3.,3.) circle (2.5pt);
		\draw [fill=ffffff] (3.,2.) circle (2.5pt);
		\draw [fill=ffffff] (6.,3.) circle (2.5pt);
		\draw [fill=ffffff] (6.,2.) circle (2.5pt);
		\draw [fill=ffffff] (4.,1.) circle (2.5pt);
		\draw [fill=ffffff] (5.,3.) circle (2.5pt);
		\draw [fill=ffffff] (5.,2.) circle (2.5pt);
		\draw [fill=ffffff] (10.,3.) circle (2.5pt);
		\draw [fill=ffffff] (10.,2.) circle (2.5pt);
		\draw [fill=ffffff] (11.,3.) circle (2.5pt);
		\draw [fill=ffffff] (11.,2.) circle (2.5pt);
		\draw [fill=ffffff] (13.,3.) circle (2.5pt);
		\draw [fill=ffffff] (13.,2.) circle (2.5pt);
		\draw [fill=ffffff] (14.,3.) circle (2.5pt);
		\draw [fill=ffffff] (14.,2.) circle (2.5pt);
		\draw [fill=ffffff] (11.54,1.) circle (2.5pt);
		\draw [fill=ffffff] (12.54,1.) circle (2.5pt);
		\end{scriptsize}
		\end{tikzpicture}
	\end{center}
\vspace{-3em}
	\caption{Graphs $S_4^1$ and $S_4^2$.} 
	\label{fig:K'_{1,n}}
\end{figure}
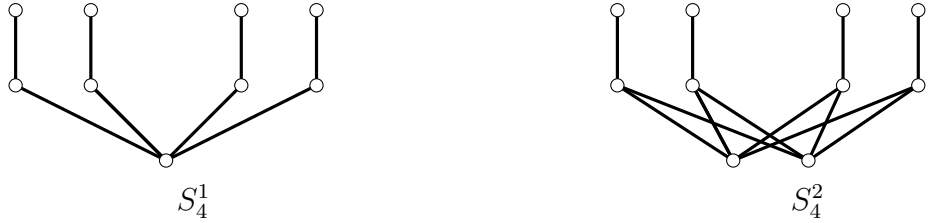

Concerning the largest possible finite values of the SMBD-numbers, we have the following:
\begin{proposition} \label{prop:n/2}
	Let $G$ be a graph on $n$ vertices.
	\begin{itemize}
		\item[$(i)$] If $\gsmb'(G) <\infty$, then $\gsmb'(G)\le \left\lceil \frac{n}{2}\right\rceil$.
		\item[$(ii)$] If $\gsmb(G) <\infty$, then $\gsmb(G)\le \left\lfloor \frac{n}{2}\right\rfloor.$
	\end{itemize}
	Moreover, both upper bounds are sharp.
\end{proposition}
\proof  Staller and Dominator alternately  select vertices in the  MBD game. It readily follows that if Staller can win in an S-game or D-game, then she plays at most $\lceil n/2 \rceil$ or at most $\lfloor n/2 \rfloor $ vertices, respectively.

To show the sharpness of the upper bounds, first consider the S-game on a subdivided star $S_k^1$ which obtained from a star $S_k$ by subdividing each of its $k$ edges exactly once. Let Dominator play according the following strategy. If Staller plays a support vertex, Dominator replies by playing the adjacent leaf, and vice versa. If Staller plays the center of the star, Dominator plays one of its neighbors (if possible). This strategy of Dominator ensures that Staller cannot claim the closed neighborhood of any support vertices or leaves. As claiming the closed neighborhood of the center needs at least $k+1$ moves, we have  $\gsmb'(S_k^1)\ge k+1$. On the other hand, if Staller plays support vertices, it forces Dominator to respond with playing the adjacent leaves and Staller can win the game by playing the last leaf or the center of the star. Therefore $\gsmb'(S_k^1)<\infty$, and by the upper bound $\lceil \frac{n}{2} \rceil= k+1$ we conclude $\gsmb(S_k^1)=k+1$. This proves that for every odd integer $n$, there exists a graph with $\gsmb'= \lceil n/2 \rceil$. If $n$ is even, let $k= (n-2)/2$. Supplementing $S_k^1$ with a new vertex that is adjacent to all support vertices, we obtain a graph  $S_k^2$ on $n=2k+2$ vertices which satisfies $\gsmb'(S_k^2)=k+1$. To obtain sharp examples for the D-game, we may add an isolated vertex to the graphs $S_k^1$ and $S_k^2$.
\qed

Even if the upper bounds in Proposition~\ref{prop:n/2} are tight, we do not know sharp examples with $\delta(G) \ge 2$. We conjecture that the inequalities can be improved as follows.
\begin{conjecture}
	Let $G$ be a graph $n$ vertices.
\begin{itemize}
	\item[$(i)$] If $\gsmb'(G) <\infty$, then $\gsmb'(G)\le \left\lceil \frac{n}{2}\right\rceil - \delta(G) +1$.
	\item[$(ii)$] If $\gsmb(G) <\infty$, then $\gsmb(G)\le \left\lfloor \frac{n}{2}\right\rfloor - \delta(G)+1.$
\end{itemize}
\end{conjecture}

\section*{Acknowledgements}
We are grateful to Sandi Klav\v zar for helpful discussions on the topic of the manuscript. The first author acknowledges the financial support from the Slovenian Research Agency under the project N1-0108.

\end{document}